\documentclass[10pt]{amsart}
\usepackage{appendix}
\usepackage{graphicx,enumitem,dsfont}
\usepackage[usenames]{color}
\usepackage[colorlinks]{hyperref}
\usepackage{subfigure}

\newcommand{\Z}{{\mathbb{Z}} }

\newcommand{\R}{\mathbb{R}}

\newcommand{\U}{{\mathbf{U}} }
\newcommand{\V}{{\mathbf{V}} }

\newcommand{\Comment}[1]{}

\renewcommand{\i}{\ifmmode\mathit{\mathchar"7010 }\else\char"10 \fi}
\renewcommand{\j}{\ifmmode\mathit{\mathchar"7011 }\else\char"11 \fi}

\newcommand{\Dt}{\Delta t}

\newcommand{\norm}[1]{\left\|#1\right\|}
\newcommand{\abs}[1]{\left|#1\right|}

\newtheorem{theorem}{Theorem}[section]

\newtheorem{proposition}{Proposition}[section]
\newtheorem{lemma}{Lemma}[section]

\newtheorem{remark}{Remark}[section]
\theoremstyle{definition} 

\newtheorem*{maintheorem*}{Main Theorem}

\allowdisplaybreaks

\numberwithin{equation}{section}
\numberwithin{figure}{section}
\numberwithin{table}{section}

\newcounter{asnr}

{\ifnum\value{asnr}=0 \stepcounter{asnr} 
  \begin{enumerate}[label=\textbf{A}.\arabic{enumi}]
    \else
    \begin{enumerate}[label=\textbf{A}.\arabic{enumi},resume] \fi}
{\end{enumerate}}

\title[Kawahara equation]{Error estimate for a fully discrete spectral scheme for Korteweg-de Vries-Kawahara equation.}

\author[U. Koley]{U. Koley} \address[Ujjwal Koley]{\newline Centre of
  Mathematics for Applications (CMA) \newline University of
  Oslo\newline P.O. Box 1053, Blindern\newline N--0316 Oslo, Norway}
\email[]{toujjwal@gmail.com}

\date{\today}

\begin{document}

\begin{abstract}
  We are concerned with the convergence of spectral method for the numerical solution of the initial-boundary value problem associated to the Korteweg-de Vries-Kawahara equation (in short Kawahara equation), which is a transport equation perturbed by dispersive terms of $3$rd and $5$th order. This equation appears in several fluid dynamics problems. It describes the evolution of small but finite amplitude long waves in various problems in fluid dynamics. These equations are discretized in space by the standard Fourier-Galerkin spectral method and in time by the explicit leap-frog scheme. For the resulting fully discrete, conditionally stable scheme we prove an $L^2$-error bound of spectral accuracy in space and of second-order accuracy in time.
\end{abstract}

\maketitle

\section {Introduction}
In this paper, we analyze the numerical approximation by Fourier spectral methods to the
Korteweg-de Vries-Kawahara (briefly Kawahara) equation with periodic solutions:
\begin{equation}
\begin{cases}
\label{eq:kawahara}
u_t  = - u u_x - u_{xxx} + u_{xxxxx}, \quad & (x,t)\in \R \times [0,\infty),\\
u(x,t) = u(x+ 2\pi,t), \quad & (x,t)\in \R \times [0,\infty), \\
u (x,0) = f(x),\quad & x \in \R,
\end{cases}
\end{equation}
where the initial condition $f$ is a given real valued $2 \pi$-periodic function.

It is well known that the one-dimensional waves of small but finite amplitude in dispersive systems (e.g., the magneto-acoustic waves in plasmas, the shallow water waves, the lattice waves and so on) can be described by the Korteweg-de Vries (KdV in short) equation, given by
\begin{equation}
\begin{aligned}
\label{eq:kdv}
u_t  = - u u_x - u_{xxx} , 
\end{aligned}
\end{equation}
which admits either compressive or rarefactive steady solitary wave solution (by a solitary water wave, we mean a travelling wave solution of the water wave equations for which the free surface approaches a constant height as $|x| \rightarrow \infty$) according to the sign of the dispersion term (the third order derivative term). Under certain circumstances, however, it might happen that the coefficient of the third order derivative in the KdV equation becomes small or even zero. In that case one has to take account of the higher order effect of dispersion in order to balance the nonlinear effect. In such cases one may obtain a generalized nonlinear dispersive equation, known as Kawahara equation, which has a form of the KdV equation with an additional fifth order derivative term given by \eqref{eq:kawahara}.
The Kawahara equation is an important nonlinear dispersive equation. It describes solitary wave propagation in media in which the first-order dispersion is anomalously small. A more specific physical background of this equation was introduced by Hunter and Scheurle \cite{hun}, where they used it to describe the evolution of solitary waves in fluids in which the Bond number is less than but close to $\frac{1}{3}$ and the Froude number is close to $1$. In the literature this equation is also referred to as the fifth order KdV equation or singularly perturbed KdV equation. The fifth order term $\partial_x^5 u$ is called the Kawahara term.
There has been a great deal of work on solitary wave solutions of the Kawahara equation \cite{kawa,kic,kato1,pon,gleb} over the past thirty years. It is found that, similarly to the KdV equation, the Kawahara equation also has solitary wave solutions which decay rapidly to zero as $t\to \infty$, but unlike the KdV equation whose solitary wave solutions are non-oscillating, the solitary wave solutions of the Kawahara equation have oscillatory trails. This shows that the Kawahara equation is not only similar but also different from the KdV equation in the properties of solutions, like what happens between the formulations of this equation and the KdV equation. The strong physical background of the Kawahara equation and such similarities and differences between it and the KdV equation in both the form and the behavior of the solution render the mathematical treatment of this equation particularly interesting. The Cauchy problem for Kawahara equation has been studied by a few authors \cite{bia,ver,kato2,hua,shang}. It has been shown that the Cauchy problem has a local solution $u \in C([-T,T];H^r(\R))$ if $f \in H^r(\R)$ and $r>-1$. This local result combined with the energy conservation law yields that \eqref{eq:kawahara} has a global solution $u \in C([-\infty,\infty];L^2(\R))$ if $f \in L^2(\R)$. Well-posedness results can be found in \cite{bia}.

Being integrable, Kawahara equation \eqref{eq:kawahara} has infinitely many invariants.
Below we will state only first three of them.
\begin{lemma}
\label{lemma:invariants}
There exists a unique solution to \eqref{eq:kawahara}.
Moreover this solution conserves the first three energy integrals, namely
\begin{align}
&(\partial / \partial t) \Biggl[ \int_{0}^{2\pi} u(x,t)\,dx \Biggr] = 0,\label{eq:inv_1} \\
&(\partial / \partial t) \Biggl[ \int_{0}^{2\pi} u^2(x,t)\,dx \Biggr ] = 0, \label{eq:inv_2} \\
&(\partial / \partial t) \Biggl [\int_{0}^{2\pi} \left( \frac{1}{3} u^3 -  u_x^2 - u_{xx}^2 \right)\,dx \Biggr] = 0.\label{eq:inv_3}
\end{align}
\end{lemma}

\begin{proof} \ The invariance of these expressions can be shown for smooth solutions by using
periodicity. For the sake of completeness, we will give a proof.

In order to show \eqref{eq:inv_1}, let us integrate \eqref{eq:kawahara} in space. We get
\begin{align*}
(\partial / \partial t) \int_{0}^{2\pi} u \,dx  + (1/2) \int_{0}^{2\pi} (u^2)_x \,dx 
+ \int_{0}^{2\pi} u_{xxx} \,dx  -\int_{0}^{2\pi} u_{xxxxx} \,dx = 0,
\end{align*}
using the periodicity of $u$, we deduce then \eqref{eq:inv_1}.

To prove \eqref{eq:inv_2} we start by multyplying the equation \eqref{eq:kawahara} by $u $ and integrate by parts in space, yields
\begin{align*}
\int_{0}^{2\pi} uu_t \,dx & = \int_{0}^{2\pi} -u^2 u_x - uu_{xxx} +  uu_{xxxxx} \,dx \\
& = -\int_{0}^{2\pi} (\frac{1}{3} u^3)_x \,dx - \int_{0}^{2\pi} (uu_{xx} - \frac{1}{2}u_x^2)_x \,dx - \int_{0}^{2\pi} u_x u_{xxxx}\,dx \\
&= -\int_{0}^{2\pi} (\frac{1}{3} u^3)_x \,dx - \int_{0}^{2\pi} (uu_{xx} - \frac{1}{2}u_x^2)_x \,dx \\
& \qquad - \int_{0}^{2\pi} (u_xu_{xxx} - \frac{1}{2}u_{xx}^2)_x \,dx. 
\end{align*}
Again using the periodicity of $u$, we can establish \eqref{eq:inv_2}. 

To prove \eqref{eq:inv_3}, we start by multyplying \eqref{eq:kawahara} by $u^2 $ and integrate by parts in space, yields
\begin{align*}
\int_{0}^{2\pi} u^2u_t \,dx & = \int_{0}^{2\pi} - u^3 u_x  - u^2u_{xxx} +  u^2u_{xxxxx} \,dx \\
& = -\int_{0}^{2\pi} (\frac{1}{4} u^4)_x \,dx + 2 \int_{0}^{2\pi} (uu_x)u_{xx} \,dx -  2 \int_{0}^{2\pi} (uu_x) u_{xxxx}\,dx\\
& =2 \int_{0}^{2\pi} [ -u_t -  u_{xxx} + u_{xxxxx} ]u_{xx} \,dx \\
& \qquad \qquad - 2 \int_{0}^{2\pi}[ - u_t - u_{xxx} + u_{xxxxx} ] u_{xxxx}\,dx \\
& =   2 \int_{0}^{2\pi} u_{tx}u_x \,dx - 2 \int_{0}^{2\pi} u_{xx}u_{xxx} \,dx + 2 \int_{0}^{2\pi} u_{xx}u_{xxxxx} \,dx \\
 - 2 & \int_{0}^{2\pi} u_{tx}u_{xxx} \,dx + 2  \int_{0}^{2\pi} u_{xxx}u_{xxxx} \,dx - 2 \int_{0}^{2\pi} u_{xxxx}u_{xxxxx} \,dx \\
&= 2 \int_{0}^{2\pi} u_{tx}u_x \,dx + 2 \int_{0}^{2\pi} u_{txx}u_{xx} \,dx.
\end{align*}
From this we can conclude that
\begin{align*}
\partial / \partial t \Biggl[\int_{0}^{2\pi} \left( \frac{1}{3}u^3 -  u_x^2 - u_{xx}^2 \right) \,dx \Biggr] = 0.
\end{align*}
\end{proof}

There has been a great deal of work on the Fourier-Galerkin spectral method for 
the KdV equations \cite{bona,baker,maday}. Also, spectral methods for initial- and
periodic boundary value problems for nonlinear wave equations with nonlocal
dispersive terms has been studied by many authors \cite{pelloni}. 
In this paper, we prove error estimates for a simple spectral fully discrete scheme
that we use to approximate spatially periodic solutions of Kawahara equation. We
first discretize in space using the standard Fourier-Galerkin spectral method,
which is easily shown to preserve the first three invariants of the equation.
We prove that the solution of the semi discrete problems converges in $L^2$,
on bounded temporal intervals, to the solution of the corresponding 
continuous problem at the spectral rate $N^{-r}$; here the number of Fourier
points is $2N+1$ and $r$ is the order of the Sobolev space in which the
solution is supposed to belong for $t=0$. 

We then discretize in time the ODE system that results from the spectral 
semidiscretization using two different methods. First, a second order accurate explicit scheme,
the leap frog method and secondly, a second order accurate semi-implicit scheme,
the Crank-Nicholson method. We prove the expected $\mathcal{O}(\Dt^{2})$ error
bound in $L^2$ for both these temporal discretization, for suitably accurate
initial conditions and under the stability requirement that $\Dt N^5$ and 
$\Dt N$ respectively are
sufficiently small. The same type of mesh restriction is required for 
stability of any explicit temporal discretization of the stiff ODE
semi discrete system under consideration.

The rest of the paper is organized as follows: In section $2$, we give all the 
necessary preliminary results. In section $3$, we consider a semi discrete 
Fourier-Galerkin scheme for the initial-boundary value problem corresponding
to \eqref{eq:kawahara} and prove an error estimate. We prove that the solution of the
semi discrete problem converges in $L^2$. In section $4$, we consider a fully
discrete explicit Fourier-Galerkin scheme for the initial-boundary value problem 
corresponding to \eqref{eq:kawahara} and prove an error estimate under the stability
condition that $\Dt N^5$ is sufficiently small. Finally, in section $5$, we consider a fully
discrete semi explicit Fourier-Galerkin scheme for the initial-boundary value problem 
corresponding to \eqref{eq:kawahara} and prove an error estimate under the stability
condition that $\Dt N$ is sufficiently small.   

\section{ Notation and Preliminary results}
\label{sec:notation}

We consider functions that are periodic of period $2 \pi$. The function spaces
we use here are $L^2$ and the Sobolev spaces $H^r$ for integer $r \ge 0$.
These spaces will always be considered on $[-\pi, \pi]$ and their elements
will be periodic functions. We denote by $(.,.)$ the standard $L^2$ inner
product; this yields a norm in $L^2$ which we denote by $\norm{.}$.
The norm in $H^r$, denoted by $\norm{.}_r$, is defined by
\begin{align}
\norm{f}_r = \left( \sum_{k \in \Z} (1 + k^2)^r \abs {\hat{f}(k)}^2\right)^{1/2}.
\end{align}
By $\norm{.}_{\infty}$, we denote the norm of $L^{\infty} = L^{\infty} [-\pi, \pi]$.

As usual, we denote by $\hat{f}(k)$, $k \in \Z$, the Fourier coefficients of $f$:
\begin{align*}
\hat{f}(k) = \frac{1}{2\pi} \int_{-\pi}^{\pi} e^{-ikx} f(x) \,dx.
\end{align*}

We recall that the Fourier coefficients of the pointwise product are given
by the convolution of the Fourier coefficients of $f$ and $g$, defined by
\begin{align*}
(\hat{f} * \hat{g})(k) = \sum_{m,n \in \Z; m+n=k} \hat{f}(m)\hat{g}(n). 
\end{align*}

We also need to consider discrete analogues of the quantities defined above.
To this end, for $N$ a positive integer, consider the space $S_N$ defined by
\begin{align}
S_N = \text{span} \lbrace { exp(ikx): k \in \Z, -N \le k \le N} \rbrace.
\end{align}

Let $P_N$ denote the $L^2$ orthogonal projection onto $S_N$. The projection has
the following approximation properties, whose proof is standard.

\begin{proposition}
Given integers $ 0 \le s \le r $, there exists a constant $C$ independent of $N$ such that,
for any $f \in H^r$,
\begin{align}
&\norm { f- P_Nf }_{s} \le C N^{s-r} \norm{f}_{r},\label{eq:pro_1} \\
&\norm { f- P_Nf }_{\infty} \le C N^{1/2-r} \norm{f}_{r}, \quad r \ge 1. \label{eq:pro_2}
\end{align}
\end{proposition}
\begin{proposition}
Given integers $ 0 \le s \le r $, there exists a constant $C$ independent of $N$ such that,
for any $\psi \in S_N$,
\begin{align}
\norm {\psi}_r \le C N^{r-s} \norm{\psi}_s, \quad \text{and} \quad \norm{\psi}_{\infty} \le C N^{1/2} \norm{\psi}.\label{eq:pro_3}
\end{align}
\end{proposition}

We use many times in the proofs a Sobolev-type inequality, which states that
there exists a constant $C>0$ such that for all $ f \in H^1$
\begin{align}
&\norm {f}_{\infty} \le C  \norm{f}^{1/2} \norm{f}_{1}^{1/2}.\label{eq:pro_3}
\end{align}

As we have already mentioned in the introduction, we shall consider the initial- and
periodic boundary-value problem for the Kawahara equation: We seek a real-valued
function $u(x,t)$, $2 \pi$ periodic in $x$ and satisfying
\begin{equation}
\begin{cases}
\label{eq:main}
u_t + u u_x + u_{xxx} = u_{xxxxx}, \quad & x \in [-\pi, \pi], \quad t \ge 0,\\
u(x,0) = f(x), \quad & x \in [-\pi, \pi].
\end{cases}
\end{equation}
Here $f(x)$ is a real-valued, $2 \pi$ periodic function.

We now going to state a existence and uniqueness result for the solution
of \eqref{eq:main}. See \cite{ujjwal} for a proof.

\begin{theorem}
Let $f$ be in $H^r$, with $r \ge 5$. Then there exists a unique solution $u$ of
\eqref{eq:main} in $H^r$, which belongs to the space $ C^k( [0, \infty); H^{r-5k})$ with $r -5k \ge -1$, i.e., is such that its temporal derivatives up to order $k$ exist and are continuous and bounded on $[0, \infty)$ with values in $H^{r -5k}$.
\end{theorem}

\section{Semi-discrete approximation:}
 
In this section we analyze a Fourier-Galerkin scheme for the discretization of \eqref{eq:kawahara} in the spatial
variable.

The semi discrete Fourier-Galerkin (spectral) approximation to \eqref{eq:kawahara} is a map $\U$ from $[0, \infty)$ to the real-valued elements of $S_N$ such that, for all $\phi \in S_N$:
\begin{equation}
\label{eq:kawa_semi}
\begin{cases}
\left( \U_t + \U\U_x +\U_{xxx} - \U_{xxxxx}, \phi \right)=0, \quad t\ge 0, \\
\U(0)= P_Nf,
\end{cases}
\end{equation} 
where $P_N$ is the orthogonal projection of $L^2$ onto $S_N$.

By choosing $\phi = e^{ikx}$ for $k = -N,\cdots,N$, we see that \eqref{eq:kawa_semi} is 
an initial-value problem for an ODE system for Fourier coefficients $\hat{\U}(k,t)$ of
$\U$. Since $\U$ is real, these coefficients must satisfy the condition 
$\hat{\U}(k,t) = \overline{\hat{\U}(k,t)}$ and the equation
\begin{equation}
\label{eq:kawa_fourier}
\begin{aligned}
\hat{\U}_t(k,t) & = \frac{-ik}{2} \hat{\U}* \hat{\U}(k,t) - k^3\hat{\U}(k,t) - k^5\hat{\U}(k,t) , \quad -N \le k \le N\\
\hat{\U}(k,0) & = \hat{f}(k). 
\end{aligned}
\end{equation}

The right hand side of the system \eqref{eq:kawa_fourier} is Lipschitz continuous, at least
locally, with respect to $l^2$ norm. Hence, the existence of a maximal time $t_0$, $0 < t_0 \le T$
such that, for all $t < t_0$, there exists a unique solution $U(t)$ to problem \eqref{eq:kawa_semi}
is a classical result of the theory of differential systems. The problem is to get the existence
for an arbitrary time $t_0$, or equivalently to prove that one can take $t_0 = T$. This result
is a consequence of the fact that \eqref{eq:kawa_semi} is conservative in $L^2$, which ensures
that the solution cannot blow-up.

Now we will present the main properties enjoyed by the Fourier-Galerkin approximation
\eqref{eq:kawa_semi}. More precisely, this semidiscretization preserves the discrete
analogues of the first three invariants of \eqref{eq:main}.

\begin{lemma}
\label{lemma:inv_discrete}
There exists a unique solution $\U$ to problem \eqref{eq:kawa_semi}.
Moreover this solution conserves the first three energy integrals of Kawahara equation, namely
\begin{align}
&(\partial / \partial t) \Biggl[ \int_{-\pi}^{\pi} \U(x,t)\,dx \Biggr] = 0,\label{eq:inv_1d} \\
&(\partial / \partial t) \Biggl[ \int_{-\pi}^{\pi} \U^2(x,t)\,dx \Biggr ] = 0, \label{eq:inv_2d} \\
&(\partial / \partial t) \Biggl [\int_{-\pi}^{\pi} \left( \frac{1}{3} \U^3 -  \U_x^2 - \U_{xx}^2 \right)(x,t)\,dx \Biggr] = 0.\label{eq:inv_3d}
\end{align}
\end{lemma}

\begin{proof}
First of all, we have already discussed about the existence of a unique solution to \eqref{eq:kawa_semi}.
Now in order to show \eqref{eq:inv_1d}, let us first choose $\phi =1$ as a test function in 
\eqref{eq:kawa_semi}. We get
\begin{align*}
(\partial / \partial t) \int_{-\pi}^{\pi} \U(x,t) \,dx + (1/2) \int_{-\pi}^{\pi} \U^2_x(x,t) \,dx
& + \int_{-\pi}^{\pi} \U_{xxx}(x,t) \,dx \\
& - \int_{-\pi}^{\pi} \U_{xxxxx}(x,t) \,dx = 0,
\end{align*}
using the periodicity of $\U$, we deduce then \eqref{eq:inv_1d}.
To prove \eqref{eq:inv_2d}, we choose $\phi = \U$ in \eqref{eq:kawa_semi}. We obtain
\begin{align*}
(\partial / \partial t) \int_{-\pi}^{\pi} \U^2(x,t) \,dx + (1/3) \int_{-\pi}^{\pi} \U^3_x(x,t) \,dx
& + \int_{-\pi}^{\pi} (\U \U_{xxx})(x,t) \,dx \\
& - \int_{-\pi}^{\pi} (\U \U_{xxxxx})(x,t) \,dx = 0.
\end{align*}
Integrating by parts and using the periodicity of $\U$ yields
\begin{align*}
\int_{-\pi}^{\pi} (\U \U_{xxx})(x,t) \,dx = - (1/2) \int_{-\pi}^{\pi} (\U_{x}^2)_x(x,t) \,dx = 0.
\end{align*}
Similarly, we have
\begin{align*}
\int_{-\pi}^{\pi} (\U \U_{xxxxx})(x,t) \,dx  = 0, \quad \text{and} \quad \int_{-\pi}^{\pi} \U^3_{x}(x,t) \,dx  = 0.
\end{align*}
Hence, we deduce \eqref{eq:inv_2d}.

\noindent We derive now \eqref{eq:inv_3d} by choosing $\phi = P_N(\U^2)$ in \eqref{eq:kawa_semi}. As a result we obtain
\begin{align*}
\int_{-\pi}^{\pi} \U_t P_N \U^2 \,dx + \int_{-\pi}^{\pi} \U\U_x P_N\U^2 \,dx
& + \int_{-\pi}^{\pi} \U_{xxx} P_N\U^2 \,dx \\
& - \int_{-\pi}^{\pi} \U_{xxxxx} P_N\U^2 \,dx = 0.
\end{align*}
Now as $\U_t$ is an element of $S_N$, we have
\begin{align*}
\int_{-\pi}^{\pi} \U_t P_N \U^2 \,dx = \int_{-\pi}^{\pi} \U_t \U^2 \,dx = (1/3)(\partial / \partial t) \int_{-\pi}^{\pi} \U^3 \,dx.
\end{align*}
On the other hand,
\begin{align*}
\int_{-\pi}^{\pi} \U_{xxx} P_N\U^2 \,dx  & = \int_{-\pi}^{\pi} \U_{xxx} \U^2 \,dx = - 2\int_{-\pi}^{\pi} \U_{xx} (\U\U_x) \,dx  \\
& = - \int_{-\pi}^{\pi} \U_{xx}\left( -\U_t -\U_{xxx} + \U_{xxxxx}\right)  \,dx
= -(\partial / \partial t) \int_{-\pi}^{\pi} \U_x^2 \,dx. 
\end{align*}
Similarly we have 
\begin{align*}
 - \int_{-\pi}^{\pi} \U_{xxxxx} P_N\U^2 \,dx  = -(\partial / \partial t) \int_{-\pi}^{\pi} \U_{xx}^2 \,dx.
\end{align*}
Finally, using the fact that $P_N$ commutes with differentiation, we obtain
\begin{align*}
&\int_{-\pi}^{\pi} \U\U_x P_N\U^2 \,dx = \int_{-\pi}^{\pi} \left(\frac{1}{2}\U^2\right)_x P_N\U^2 \,dx\\
& = - \int_{-\pi}^{\pi} \U^2 P_N (\U\U_x) \,dx = - \int_{-\pi}^{\pi} \U\U_x P_N\U^2 \,dx,  
\end{align*}
consequently, we have 
\begin{align*}
\int_{-\pi}^{\pi} \U\U_x P_N\U^2 \,dx = 0.
\end{align*}
Combining all the results above we get \eqref{eq:inv_3d}.
\end{proof}

We shall now state the following theorem:

\begin{theorem}
\label{theo:main}
The semi-discrete scheme \eqref{eq:kawa_semi} has an unique solution $\U$ for $t \ge 0$. Let
$u(x,t)$ be the solution of \eqref{eq:kawahara} corresponding to the initial data $u_0 \in H^r$.
Then there exists a time $T > 0$, and a constant $C > 0$, independent of $N$, such that
\begin{equation}
\label{eq:mainerror}
\begin{aligned}
\max_{0 \le t \le T} \norm {u-\U} \le \frac{C}{N^{r}}.
\end{aligned}
\end{equation}
\end{theorem}

Before giving the proof of Theorem ~\ref{theo:main}, we define a different semi-discretization
based on a linearization of \eqref{eq:kawa_semi}. The convergence of this approximation will
serve as an intermediate step in the proof of the convergence of the original scheme.

To this end, we linearize \eqref{eq:kawa_semi} as follows: Given a solution of \eqref{eq:kawahara},
corresponding to initial data $u_0$ in $H^r$, we look for a function $V \in S_N$, which for all
$\phi \in S_N$, satisfies
\begin{equation}
\label{eq:kawa_inter}
\begin{cases}
\left( V_t + uV_x +V_{xxx} - V_{xxxxx}, \phi \right)=0, \quad t\ge 0, \\
V(0)= P_Nu_0,
\end{cases}
\end{equation}

\begin{lemma}
\label{lemma:inter}
Let $u(x,t)$ be a solution of \eqref{eq:kawahara} corresponding to initial data $u_0 \in H^r$.
Then there exists a unique solution $V$ of \eqref{eq:kawa_inter} for all $t \ge 0$.
Moreover, given $ 0 \le t \le T$, there exists a constant $C$ independent of $N$ such that
\begin{equation}
\label{eq:intererror}
\begin{aligned}
\max_{0 \le t \le T} \norm {u-V} \le \frac{C}{N^{r}}.
\end{aligned}
\end{equation}
\end{lemma}

\begin{proof}
The existence of an unique local solution of \eqref{eq:kawa_inter} is again a consequence
of standard ODE theory. To see that we have global existence, we resort to a stability
result in the $L^2$ norm. Choosing $\phi = V$ in \eqref{eq:kawa_inter}, we obtain
\begin{align*}
\frac{1}{2} \frac{d}{dt} \norm{V}^2 \le \norm{u_x}_{\infty} \norm{V}^2.
\end{align*}
Sobolev's inequality and the fact that $u_0 \in H^2$ imply that $\max_{t \ge 0} \norm{u_x}_{\infty} \le C$;
thus by the Gronwall inequality, there exists $C$ such that $\max_{0 \le t \le T} \norm{V} \le C$, and 
we can extend the local solution to a solution on every bounded interval $[0,T]$.
Note that using the same arguments and choosing $\phi = V_{xx}$, we conclude that
$\max_{0 \le t \le T} \norm{V_x} \le C$.

Now set $\rho = P_N u -V$, then $ u- V = u -P_Nu + P_N u -V = u -P_Nu + \rho$.
The hypothesis on $u$ imply that $\norm {u -P_Nu} \le C/N^{r}$. Thus we only need to 
estimate $\rho$. Observe that, $\rho$ is an element of $S_N$ satisfying the equation
\begin{align*}
\left( \rho_t + \rho_{xxx} - \rho_{xxxxx} , \phi \right) = -\left( P_N (uu_x), \phi \right) +(\rho uV_x, \phi).
\end{align*}
Now observe that,
\begin{align*}
u u_x - u V_x = (u -P_Nu)u_x + (P_Nu -V)u_x + V u_x -uV_x
\end{align*}
Then, the choice $\phi = \rho$ yields,
\begin{align*}
\frac{1}{2} \frac{d}{dt} \norm{\rho}^2 \le \norm{u_x}_{\infty} \norm{u-P_Nu}\norm{\rho} + \norm{u_x}_{\infty} \norm{\rho}^2
+ \norm{u_x}_{\infty} \norm{\rho} \norm{V} + \norm{u}_{\infty} \norm{\rho} \norm{V_x} ,
\end{align*}
so that by the arithmetic-geometric inequality we obtain,
\begin{align*}
\frac{1}{2} \frac{d}{dt} \norm{\rho}^2 \le  \frac{C}{(N^r)^2} + C \norm{\rho}^2.
\end{align*}
But since $\rho(x,0)=0$, hence we obtain by using Grownwall's inequality
\begin{align*}
\max_{0 \le t \le T} \norm {\rho} \le \frac{C}{N^{r}},
\end{align*}
and this completes the proof of the lemma.
\end{proof}

\begin{proof}{(of Theorem ~\ref{theo:main})}\\
We have already proved the existence and uniqueness of the semi-discrete 
solution $\U$.
Now set $e = V -\U$. Then $u-\U = u-V + e$.
In view of \eqref{eq:intererror}, we only need an estimate for $e$. Now 
observe that $e$ satisfies
\begin{align*}
\left( e_t + e_{xxx} - e_{xxxxx}, \phi \right)= - (uV_x -\U\U_x , \phi).
\end{align*}
Now using the fact that
\begin{align*}
uV_x -\U\U_x = (u-V)V_x + (eV)_x - ee_x
\end{align*}
and choosing $\phi = e$, observing that $(e,ee_x)=0$, yields,
\begin{equation}
\label{eq:est1}
\begin{aligned}
\frac{1}{2} \frac{d}{dt} \norm{e}^2 \le \norm{V_x}_{\infty} \norm{u-V} \norm{e} + \frac{1}{2} \norm{V_x}_{\infty}\norm{e}^2.
\end{aligned}
\end{equation}
We now use the following {\bf inverse} inequalities: for $0 \le s \le r$ and $\psi \in S_N$
\begin{align*}
\norm {\psi}_r \le C N^{r-s} \norm{\psi}_s, \qquad \norm{\psi}_{\infty} \le C N^{1/2} \norm{\psi}.
\end{align*} 
Since $(P_Nu -V)$ is in $S_N$,
\begin{align*}
\norm{(P_Nu -V)_x}_{\infty} \le C \norm{P_Nu -V}_{1}^{1/2} \norm{P_Nu -V}_{2}^{1/2} \le C N^{3/2} \norm{P_Nu -V}.
\end{align*}

\noindent Consequently, we have 
\begin{align*}
\norm{P_Nu -V}_{\infty} \le C N^{1/2 -r}, \\
\norm{(P_Nu -V)_x}_{\infty} \le C N^{3/2 -r}.
\end{align*}
Using Sobolev's inequality along with the approximation properties of the projection $P_N$, we have
\begin{align*}
\norm{V_x}_{\infty} & \le \norm{u_x}_{\infty} + \norm{(u -P_Nu)_x}_{\infty} + \norm{(P_Nu -V)_x}_{\infty}\\
& \le C + C \norm{u -P_Nu}_{2} + \frac{C}{N^{r-3/2}}\\
& \le C + \frac{C}{N^{r-2}} + \frac{C}{N^{r-3/2}}.
\end{align*}

\noindent Using $r \ge 2$, these inequalities yield
\begin{align*}
\max_{0 \le t \le T}  \norm{V_x}_{\infty} \le C,
\end{align*}
consequently (using arithmetic-geometric inequality), from \eqref{eq:est1}
\begin{align*}
\frac{1}{2} \frac{d}{dt} \norm{e}^2 \le C \left( \norm{e}^2 +  \norm{u-V}^2 \right).
\end{align*}

\noindent Finally, writing $u -V = u -P_Nu + \rho$, we have
\begin{align*}
\frac{1}{2} \frac{d}{dt} \norm{e}^2 \le C \norm{e}^2 +  \left(\frac{C}{N^{r}}\right)^2.
\end{align*}
Since $e(0)=0$, Gronwall's inequality gives,
\begin{align*}
\max_{0 \le t \le T} \norm{e} \le \frac{C}{N^{r}},
\end{align*}
which in view of \eqref{eq:intererror}, yields \eqref{eq:mainerror}.
\end{proof}

\section{Fully-discrete scheme:}

To define the fully discrete scheme, given $0 < T < \infty$, choose a time step
$\Dt$, and an integer $M$, such that $M\Dt = T$. Then for $m=0,\cdots,M$, denote 
$t_m = m \Dt$. The fully discrete solution is defined as the sequence ${\lbrace \U^m \rbrace}$
of elements of $S_N$ satisfying, for all $\phi \in S_N$ and for $m = 1,2,\cdots, M$,
the equation
\begin{equation}
\label{eq:kawa_fully}
\begin{aligned}
\left( \U^{m+1} - \U^{m-1}, \phi \right) + 2 \Dt \left( \U^m\U^m_x +\U^m_{xxx} - \U^m_{xxxxx}, \phi \right)=0.
\end{aligned}
\end{equation} 
For each $m$, $\U^m$ is an approximation of $\U(t_m)$, the semi-discrete solution $\U$
evaluated at time $t =t_m$. We also suppose that initial values $\U^0, \U^1$ have been
given in $S_N$.

\begin{theorem}
\label{theo:main_fully}
Let $\U(t)$ be the solution of the semi-discrete problem \eqref{eq:kawa_semi} and ${\lbrace \U^m \rbrace}$ be 
the solution of \eqref{eq:kawa_fully}. Suppose that $\U^0 = \U(0)$ and that $\U^1$ is computed in
such a way that
\begin{align}
\label{eq:initial}
\norm {\U^1 - \U(\Dt)} \le C \Dt^2.
\end{align}
Assume that $u_0$ is in $H^r$ with $ r \ge 16$. Then, there exists a constant
$C_1$ independent of $N$ and $\Dt$, such that if 
\begin{align}
\label{eq:cfl}
N^5 \Dt \le C_1,
\end{align}
there holds
\begin{equation}
\label{eq:fullyerror}
\begin{aligned}
\max_{0 \le m \le M} \norm {\U^m -\U(t_m)} \le C \Dt^2
\end{aligned}
\end{equation}
\end{theorem}

\begin{proof}

We see that ${\lbrace \U^m \rbrace}$ satisfies for all $\phi \in S_N$:
\begin{align}
\label{eq:eq-1}
\left( \U^{m+1} - \U^{m-1}, \phi \right) + 2 \Dt \left( \U^m\U^m_x +\U^m_{xxx} - \U^m_{xxxxx}, \phi \right)=0.
\end{align}
On the other hand, the semi-discrete solution $\U$ satisfies, for all $\phi \in S_N$:
\begin{align}
\label{eq:eq-2}
\left( \U(t_{m+1}) - \U(t_{m-1}), \phi \right) & + 2 \Dt \left( \U(t_m)\U_x(t_m) +\U_{xxx}(t_m) - \U_{xxxxx}(t_m), \phi \right) \\
& \qquad \qquad \qquad \qquad \qquad \qquad \qquad \qquad  =(\theta^m, \phi),
\end{align}
where $\theta^m$ is an element of $S_N$ given by:
\begin{align*}
\theta^m = \U(t_{m+1}) - \U(t_{m-1}) -2 \Dt \U_t(t_m).
\end{align*}
From Taylor's expansion, we have
\begin{align*}
\norm{\theta^m} \le C \Dt^3 \max_{t_{m-1} \le s \le t_{m+1}} \norm{\frac{\partial^3 \U(s)}{\partial t^3}}.
\end{align*}
Now let us define $e^m \in S_N$ as
\begin{align*}
e^m = \U^m - \U(t_m).
\end{align*}
Then $e^m $ satisfies, for all $\phi \in S_N$
\begin{align}
\label{eq:eq-3}
\left( e^{m+1} - e^{m-1}, \phi \right) + 2 \Dt \left( \U^m\U^m_x  - \U(t_m) \U_x(t_m) + e^m_{xxx} - e^m_{xxxxx}, \phi \right)=-(\theta^m, \phi).
\end{align}
Choosing $\phi =  e^{m+1} + e^{m-1}$, and adding $\norm{e^m}^2$ to both sides of \eqref{eq:eq-3}, we obtain
\begin{equation}
\begin{aligned}
\label{eq:eq-4}
& \norm{e^{m+1}}^2  + \norm{e^{m}}^2 - 2\Dt (e_{xx}^m , e_x^{m+1}) - 2\Dt (e_{xxx}^m, e_{xx}^{m+1})  \\
& = \norm{e^{m}}^2  + \norm{e^{m-1}}^2 - 2\Dt (e_{xx}^{m-1} , e_x^{m}) - 2\Dt (e_{xxx}^{m-1}, e_{xx}^m) - (\theta^m,  e^{m+1} + e^{m-1}) \\
& - 2 \Dt \left (\U_x^m e^m,  e^{m+1} + e^{m-1} \right) -2 \Dt \left (\U(t_m) e_x^m,  e^{m+1} + e^{m-1} \right),
\end{aligned}
\end{equation}
where we have used the following identity
\begin{align*}
\U^m \U_x^m - \U(t_m) \U_x(t_m) = \U_x^m e^m + \U(t_m) e_x^m.
\end{align*}
Now let us define $A^{m+1}$ by
\begin{align}
\label{eq:eq-5}
A^{m+1} = \norm{e^{m+1}}^2  + \norm{e^{m}}^2 & - 2\Dt (e_{xx}^m , e_x^{m+1}) - 2\Dt (e_{xxx}^m, e_{xx}^{m+1}) \\
& -2\Dt(\U(t_m)e_{x}^{m+1}, e^{m}).
\end{align}
Then we can rewrite \eqref{eq:eq-4} as
\begin{equation}
\begin{aligned}
\label{eq:eq-6}
A^{m+1} = A^m &- 2 \Dt \left (\U_x^m e^m,  e^{m+1} + e^{m-1} \right) +2 \Dt \left (\U_x(t_m) e^m,  e^{m+1} \right)\\
& - 2\Dt \left( (\U(t_m) -\U(t_{m-1})) e_x^m, e^{m-1} \right)- (\theta^m,  e^{m+1} + e^{m-1}).
\end{aligned}
\end{equation}
In general, differentiating \eqref{eq:kawa_semi} with respect to $t$ and using the properties of $P_N$, it is straightforward 
to prove that there exist constants $\alpha_{k,s}$, independent of $N$, such that, if $r \ge s + 5k+ 1$ 
\begin{align*}
\max_{0 \le t \le T} \norm{\partial^k_t \U(t)}_s \le \alpha_{k,s}.
\end{align*}
In particular, since we assume $r \ge 16$, we have by Sobolev's inequality
\begin{align*}
\norm{\U(t_m) - \U(t_{m-1})} \le \Dt \max_{t_{m-1} \le s \le t_m} \norm{\U_t(s)}_{\infty} \le C \Dt,
\end{align*}
and 
\begin{align*}
\norm{\theta^m} \le C \Dt^3.
\end{align*}
Again, \eqref{eq:eq-6} gives after some manipulations:
\begin{equation}
\begin{aligned}
\label{eq:eq-7}
A^{m+1} \le  A^m & + C \Dt^2 N \norm{e^m} \norm{e^{m+1} + e^{m-1})} + C \Dt \norm{\U_x^m}_{\infty}\norm{e^m} \norm{e^{m+1} + e^{m-1})} \\
& + C \Dt \norm{e^m} \norm{e^{m+1}} + C \Dt^3 \norm{e^{m+1} + e^{m-1})}.
\end{aligned}
\end{equation}
Now as an internl ``inductive'' hypothesis, we assume that there exists a constant $B$, independent of $N$, such that for all $n \le m$
\begin{align*}
\norm {\U_x^m}_{\infty} \le B,
\end{align*}
then using the Cauchy-Schwartz inequality, we have
\begin{align}
\label{eq:eq-8}
A^{m+1} \le  A^m + C \Dt^5 + C \Dt ( 1 + B + N \Dt ) \left ( \norm{e^{m-1}}^2 + 2 \norm{e^m}^2 + \norm{e^{m+1}}^2 \right). 
\end{align}
On the other hand, one can show that under a stability assumption \eqref{eq:cfl}, $A^{m+1}$ is positive and 
comparable to $\norm{e^m}^2 + \norm{e^{m+1}}^2$.
Infact, 
\begin{align*}
\Dt | (e_{xx}^m , e_x^{m+1}) + (e_{xxx}^m, e_{xx}^{m+1}) & + (\U(t_m)e_{x}^{m+1}, e^{m})|\\
& \le C \Dt (N^3 + N^5) \left( \norm{e^m}^2 + \norm{e^{m+1}}^2 \right).
\end{align*}
Hence, if 
\begin{align*}
C \Dt (N^3 + N^5) \le \frac{1}{2},
\end{align*}
i.e., under a condition of the form \eqref{eq:cfl}, we have
\begin{align*}
\frac{1}{2} \left( \norm{e^m}^2 + \norm{e^{m+1}}^2 \right) \le A^{m+1} \le 2 \left( \norm{e^m}^2 + \norm{e^{m+1}}^2 \right). 
\end{align*}
Consequently by \eqref{eq:eq-8}, we have
\begin{align*}
A^{m+1} \le  A^m + C \Dt^5 + C \Dt ( C + B ) ( A^m + A^{m+1}), 
\end{align*}
which implies, for $C^* = C( C + B)$
\begin{align*}
(1 - C^* \Dt) A^{m+1} \le ( 1 + C^* \Dt) A^m + C \Dt^5 .
\end{align*}
Hence, if $\Dt$ is chosen small enough, and using the fact that 
$\norm{e^1} = \mathcal{O}(\Dt^2)$ by \eqref{eq:initial}, we obtain in the standard way
\begin{align*}
\max_{ 1 \le n \le m+1} \norm{e^n} \le C \Dt^2 e^{CT}.
\end{align*}
Observe that, the above estimate allows us to complete the inductive step. Indeed, we see that
\begin{align*}
\norm{\U_x^{m+1}}_{\infty} \le \norm{\U_x(t_{m+1})}_{\infty} + \norm{e_x^{m+1}}_{\infty} \le \norm{\U_x(t_{m+1})}_{\infty}
+ C \Dt^2 N^{3/2} e^{CT} \le B
\end{align*}
consequently, by taking $\Dt$ sufficiently small and using a condition of type \eqref{eq:cfl}, we can justify the assumption that $ \norm{\U_x^m}_{\infty} \le B$ for all $m$.
\end{proof}

\begin{remark}
In conclusion, combining the results of Theorems ~\ref{theo:main} and ~\ref{theo:main_fully}, we see that
under the hypotheses of Theorem ~\ref{theo:main_fully}, the fully discrete scheme \eqref{eq:kawa_fully}
satisfies the error estimate 
\begin{align*}
\max_{ 0 \le m \le M } \norm{ u(t_m) - \U^m} \le C \left(  \frac{1}{N^{r}} + \Dt^2 \right ).
\end{align*}
\end{remark}

\section{Semi-implicit scheme}

To define the semi-implicit fully discrete scheme, given $0 < T < \infty$, choose a time step
$\Dt$, and an integer $M$, such that $M\Dt = T$. Then for $m=0,\cdots,M$, denote 
$t_m = m \Dt$. The fully discrete solution is defined as the sequence ${\lbrace \U^m \rbrace}$
of elements of $S_N$ satisfying, for all $\phi \in S_N$ and for $m = 1,2,\cdots, M$,
the equation
\begin{equation}
\label{eq:kawa_semiimplicit}
\begin{aligned}
\left( \U^{m+1} - \U^{m}, \phi \right) + \Dt \left( \U^{m +1/2}\U^{m +1/2}_x +\U^{m +1/2}_{xxx} - \U^{m +1/2}_{xxxxx}, \phi \right)=0.
\end{aligned}
\end{equation} 
For each $m$, $\U^m$ is an approximation of $\U(t_m)$, the semi-discrete solution $\U$
evaluated at time $t =t_m$. Also, we have used the notation: $\U^{m +1/2} = \frac{1}{2}(\U^m + \U^{m+1})$.

First we shall establish the rate of convergence estimates. So, for the time
being, we assume the existence of a sequence ${\lbrace{\U^n}\rbrace}_{n=0}^{M}$
in $S_N$ satisfying \eqref{eq:kawa_semiimplicit}. Later in this section we will
discuss about the existence and uniqueness of such sequences.

\begin{theorem}
\label{theo:main_implicit}
Let $\U(t)$ be the solution of the semi-discrete problem \eqref{eq:kawa_semi} and ${\lbrace \U^m \rbrace}$ be 
the solution of \eqref{eq:kawa_fully}.
Assume that $u_0$ is in $H^r$ with $ r \ge 16$. Then, there exists a constant
$C_1$ independent of $N$ and $\Dt$, such that if 
\begin{align}
\label{eq:cfl_1}
N \Dt \le C_1,
\end{align}
there holds
\begin{equation}
\label{eq:impliciterror}
\begin{aligned}
\max_{0 \le m \le M} \norm {\U^m -\U(t_m)} \le C \Dt^2
\end{aligned}
\end{equation}
\end{theorem}

\begin{proof}
First note that, by letting $\phi = \U^{m+1/2}$ in \eqref{eq:kawa_semiimplicit}, we have
by periodicity
\begin{align*}
\frac{1}{2} \left( \norm{\U^{m+1}}^2 - \norm{\U^{m}}^2 \right) 
& = -\left( \U^{m +1/2}\U^{m +1/2}_x + \U^{m +1/2}_{xxx} - \U^{m +1/2}_{xxxxx}, \U^{m+1/2} \right) \\
&= 0,
\end{align*}
hence
\begin{align}
\label{eq:conserve}
\norm{\U^m} = \norm{\U^0}, \quad m = 0, \cdots, M.
\end{align}
We see that ${\lbrace \U^m \rbrace}$ satisfies for all $\phi \in S_N$:
\begin{align}
\label{eq:eq-11}
\left( \U^{m+1} - \U^{m}, \phi \right) + \Dt \left( \U^{m +1/2}\U^{m +1/2}_x + \U^{m +1/2}_{xxx} - \U^{m +1/2}_{xxxxx}, \phi \right)=0.
\end{align}
On the other hand, the semi-discrete solution $\U$ satisfies, for all $\phi \in S_N$:
\begin{equation}
\begin{aligned}
\label{eq:eq-21}
&\left( \U(t_{m+1}) - \U(t_{m}), \phi \right) \\
&+ \Dt \left( \U(t_{m+1/2})\U_x(t_{m+1/2}) +\U_{xxx}(t_{m+1/2}) - \U_{xxxxx}(t_{m+1/2}), \phi \right) =(\theta^m, \phi),
\end{aligned}
\end{equation}
where $\theta^m$ is an element of $S_N$ given by:
\begin{align*}
\theta^m = \U(t_{m+1}) - \U(t_{m}) -2 \Dt \U_t(t_{m+1/2}).
\end{align*}
From Taylor's expansion, we have
\begin{align*}
\norm{\theta^m} \le C \Dt^3 \max_{t_{m-1} \le s \le t_{m+1}} \norm{\frac{\partial^3 \U(s)}{\partial t^3}}.
\end{align*}
Now let us define $e^m \in S_N$ as
\begin{align*}
e^m = \U^m - \U(t_m).
\end{align*}
Then $e^m $ satisfies, for all $\phi \in S_N$
\begin{equation}
\begin{aligned}
\label{eq:eq-31}
&\left( e^{m+1} - e^{m}, \phi \right) \\
& + \Dt \left( \U^{m+1/2}\U^{m+1/2}_x  - \U(t_{m+1/2}) \U_x(t_{m+1/2}) + e^{m+1/2}_{xxx} - e^{m+1/2}_{xxxxx}, \phi \right)=-(\theta^m, \phi).
\end{aligned}
\end{equation}
Now choose $\phi = e^{m+1/2}$ in \eqref{eq:eq-31} and observe the following estimates:
\begin{align*}
\left( e^{m+1} - e^{m}, e^{m+1/2} \right) = \frac{1}{2} \left( \norm{e^{m+1}}^2 - \norm{e^{m}}^2 \right). 
\end{align*}
On the other hand, we have
\begin{align*}
&\left( \U^{m+1/2}\U^{m+1/2}_x  - \U(t_{m+1/2}) \U_x(t_{m+1/2}) , e^{m+1/2} \right) \\
& = \left( e^{m+1/2}e^{m+1/2}_x  + \U(t_{m+1/2}) e^{m+1/2}_x + \U_x(t_{m+1/2}) e^{m+1/2} , e^{m+1/2} \right),
\end{align*}
and finally, using periodicity we can conclude that
\begin{align*}
\left( e^{m+1/2}_{xxx} - e^{m+1/2}_{xxxxx}, e^{m+1/2} \right) = 0.
\end{align*}
Keeping all the above estimates in mind, we have the following:
\begin{align*}
\frac{1}{2} \left( \norm{e^{m+1}}^2 - \norm{e^{m}}^2 \right)  &\le C \Dt\norm{\U_x(t_{m+1/2})}_{\infty} \norm{e^{m+1/2}}^2 +
C \Dt^3 \norm{e^{m+1/2}}\\
& \le C \Dt^5 + C \Dt \norm{e^{m+1/2}}^2.
\end{align*}
Using a standard argument, we conclude that \eqref{eq:impliciterror} holds.
\end{proof}

We now turn into the proof of existence of a sequence ${\lbrace{\U^n}\rbrace}_{n=0}^{M}$ satisfying 
\eqref{eq:kawa_semiimplicit}. For this we shall use the following variant of the well-known
fixed point theorem of Brouwer \cite{temam}.

\begin{lemma}
\label{lemma:existence}
Let $H$ be a finite-dimensional Hilbert space with inner product $(.,.)_H$, and norm
$\norm{.}_H$. Let the map $f: H \rightarrow H$ be continuous. Suppose there exists
$\beta >0$ such that $(f(K),K)_H > 0$ for all $K$ with $\norm{K} = \beta$. Then there
exists $K^{*} \in H$, $\norm{K^{*}} \le \beta$ such that $f(K^{*})=0$.
\end{lemma}

The argument of existence of ${\lbrace{\U^n}\rbrace}_{n=0}^{M}$ proceeds in an 
inductive way. Assume that ${\lbrace{\U^j}\rbrace}_{j=0}^{n}$ exists.

For $K \in S_N$, define $f: S_N \rightarrow S_N$ by
\begin{align}
\label{eq:existence}
\left( f(K),\phi \right)= \left( K - 2\U^m,\phi \right)+ \frac{\Dt}{4}\left( KK_x,\phi \right) 
+ \frac{\Dt}{2}\left( K_{xxx} - K_{xxxxx},\phi \right), \quad \forall \phi \in S_N.
\end{align}
Such a map exists by the Riesz representation theorem; the fact that $f$ is continuous
follows easily from inverse inequalities. Furthermore, by periodicity, letting
$\phi = K$
\begin{align*}
\left( f(K),K \right)= \left( K - 2\U^m,K \right)\ge \norm{K} \left( \norm{K} - 2 \norm{\U^m}\right)
\ge \norm{K} \left( \norm{K} - 2 \norm{\U^0}\right),
\end{align*}
from \eqref{eq:conserve}. Letting $\beta > 2 \norm{\U^0}$, we deduce the existence via Lemma ~\ref{lemma:existence}
of a $K^{*} \in S_N$ such that $f(K^{*})=0$. Now letting $\U^{m+1} = K^{*} - \U^m$,
we get from \eqref{eq:existence} that
\begin{align*}
\left( \U^{m+1} - \U^{m}, \phi \right) + \Dt \left( \U^{m +1/2}\U^{m +1/2}_x +\U^{m +1/2}_{xxx} - \U^{m +1/2}_{xxxxx}, \phi \right)=0, \quad \forall \phi \in S_N,
\end{align*}
proving the existence of $\U^{m+1}$.

\noindent For uniqueness, suppose that $\V^{m+1} \in S_N$ satisfies
\begin{equation}
\label{eq:kawa_semiimplicit1}
\begin{aligned}
\left( \V^{m+1} - \V^{m}, \phi \right) + \Dt \left( \V^{m +1/2}\V^{m +1/2}_x + \V^{m +1/2}_{xxx} - \V^{m +1/2}_{xxxxx}, \phi \right)=0.
\end{aligned}
\end{equation} 
Then using $ E^{m} = \U^m - \V^m$, from \eqref{eq:kawa_semiimplicit} and \eqref{eq:kawa_semiimplicit1} we have
for all $\phi \in S_N$,
\begin{align*}
( E^{m+1} & -E^m , \phi) \\
& = - \Dt \left( \V^{m +1/2} E^{m +1/2}_x  + \U^{m +1/2}_x E^{m +1/2} + E^{m +1/2}_{xxx} - E^{m +1/2}_{xxxxx}, \phi \right).
\end{align*}
We claim that as long as $U^m$ exists, we have
\begin{align} 
\label{eq:first}
&\norm{\U^m}_{\infty} \le C \quad \text{for all}\quad m\\
&\norm{\U_x^m}_{\infty} \le C \quad \text{for all}\quad m \label{eq:second}
\end{align}
Infact, Theorem ~\ref{theo:main_implicit} and some Sobolev inequalities yields,
\begin{align*}
&\norm{U^m -u(t_m)}_{\infty} \le C N^{1/2}\norm{U^m -u(t_m)} \le C N^{1/2} (\frac{1}{N^r} + \Dt^2)\\
&\norm{U_x^m -u_x(t_m)}_{\infty} \le C N^{1/2}\norm{(U^m -u(t_m))_x} \le C N^{3/2} (\frac{1}{N^r} + \Dt^2).
\end{align*}
Using the facts $r \ge 2$, and \eqref{eq:cfl_1}, we deduce \eqref{eq:first} and \eqref{eq:second}.

\noindent Letting $\phi = E^{m+1/2}$, we deduce that
\begin{align*}
\frac{1}{2} \left( \norm{E^{m+1}}^2 - \norm{E^{m}}^2 \right)  \le C \Dt\norm{\U^{m+1/2}_x}_{\infty}& \norm{E^{m+1/2}}^2\\
& + C \Dt \norm{\V^{m+1/2}_x}_{\infty} \norm{E^{m+1/2}}^2,
\end{align*}
from the above relation, we conclude that
\begin{align*}
\norm{E^{m+1}}^2  \le \nu \norm{E^{m}}^2, \quad \text{with} \quad \nu \le 1.
\end{align*}
Now taking $\U^m = \V^m$, we see that $E^{m+1} = 0$, hence uniqueness follows.

\begin{remark}
In conclusion, combining the results of Theorems ~\ref{theo:main} and ~\ref{theo:main_implicit}, we see that
under the hypotheses of Theorem ~\ref{theo:main_implicit}, the fully discrete scheme \eqref{eq:kawa_semiimplicit}
satisfies the error estimate 
\begin{align*}
\max_{ 0 \le m \le M } \norm{ u(t_m) - \U^m} \le C \left(  \frac{1}{N^{r}} + \Dt^2 \right ).
\end{align*}
\end{remark}

\end{document}